
\magnification1200
\input amstex.tex
\documentstyle{amsppt}
\nopagenumbers
\hsize=12.5cm
\vsize=18cm
\hoffset=1cm
\voffset=2cm

\footline={\hss{\vbox to 2cm{\vfil\hbox{\rm\folio}}}\hss}

\def\DJ{\leavevmode\setbox0=\hbox{D}\kern0pt\rlap
{\kern.04em\raise.188\ht0\hbox{-}}D}

\def\txt#1{{\textstyle{#1}}}
\baselineskip=13pt
\def\hf{{\textstyle{1\over2}}}
\def\a{\alpha}\def\b{\beta}
\def\d{{\,\roman d}}
\def\e{\varepsilon}\def\E{{\roman e}}

\def\b{\beta} \def\g{\gamma}
\def\G{\Gamma}

\def\={\;=\;}

\def\zt{\zeta(\hf+it)}

\def\E{{\roman e}}

 \def\I{\Im{\roman m}\,}
\def\z{\zeta}

\def\hf{{\textstyle{1\over2}}}
\def\txt#1{{\textstyle{#1}}}

\def\le{\leqslant} \def\ge{\geqslant}
\font\tenmsb=msbm10
\font\sevenmsb=msbm7
\font\fivemsb=msbm5
\newfam\msbfam
\textfont\msbfam=\tenmsb
\scriptfont\msbfam=\sevenmsb
\scriptscriptfont\msbfam=\fivemsb
\def\Bbb#1{{\fam\msbfam #1}}

\def \CC {\Bbb C}

\font\ff=cmr8
\def\txt#1{{\textstyle{#1}}}
\baselineskip=13pt

\font\teneufm=eufm10
\font\seveneufm=eufm7
\font\fiveeufm=eufm5
\newfam\eufmfam
\textfont\eufmfam=\teneufm
\scriptfont\eufmfam=\seveneufm
\scriptscriptfont\eufmfam=\fiveeufm
\def\mathfrak#1{{\fam\eufmfam\relax#1}}

\font\tenmsb=msbm10
\font\sevenmsb=msbm7
\font\fivemsb=msbm5
\newfam\msbfam
     \textfont\msbfam=\tenmsb
      \scriptfont\msbfam=\sevenmsb
      \scriptscriptfont\msbfam=\fivemsb
\def\Bbb#1{{\fam\msbfam #1}}

\def \CC {\Bbb C}

  \def\rightheadline{{\hfil{\ff Sums of zeta squares over short intervals
 }\hfil\tenrm\folio}}

  \def\leftheadline{{\tenrm\folio\hfil{\ff
   Aleksandar Ivi\'c }\hfil}}
  \def\emptyheadline{\hfil}
  \headline{\ifnum\pageno=1 \emptyheadline\else
  \ifodd\pageno \rightheadline \else \leftheadline\fi\fi}

\topmatter

\title
ON SUMS OF SQUARES OF $|\zeta(\frac12+i\g)|$ OVER SHORT INTERVALS
\endtitle
\author   Aleksandar Ivi\'c
 \endauthor

\nopagenumbers

\medskip

\address
Aleksandar Ivi\'c, Serbian Academy of Sciences and Arts,
Knez Mihailova 35, 11000 Beograd, Serbia.
\endaddress

\bigskip

\bigskip
\email {
\tt
aleksandar.ivic\@rgf.bg.ac.rs, aivic\_2000\@yahoo.com}
\endemail
\dedicatory
\enddedicatory
\abstract
{
A
discussion involving the evaluation of the sum $$\sum_{T<\g\le
T+H}|\z(\hf+i\g)|^2$$ and some related integrals
is presented, where $\g$ denotes imaginary parts
of complex zeros of the Riemann zeta-function $\z(s)$. It is shown
unconditionally that the above sum is $\,\ll H\log^2T\log\log T\,$ for $\,T^{2/3}\log^4T \ll H \le T$.
 }
\endabstract
\endtopmatter

\document

\head
1. Introduction and statement of results
\endhead

Let $\gamma \,(>0)$ denote ordinates of complex zeros of the Riemann zeta-function $\zeta(s)$.
Consider
$$
F(T,H) := \sum_{T<\gamma\le T+H}|\zeta(\hf + i\g)|^2\qquad(1 \ll H = H(T) \le T),\leqno(1.1)
$$
so that the interval $[T, T+H]\,$ may be called ``short'' if $H = o(T)$ as $T\to\infty$.
\smallskip
A more general sum than the one in (1.1), when $H=T$, was treated by S.M. Gonek \cite{Gon}.
He proved, under the Riemann hypothesis 
(RH, that all complex zeros $\rho = \b+i\g$ of $\z(s)$ satisfy $\b = \hf$) that
$$
\sum_{0<\g\le T}\left|\z\left({1\over2}+i\left(\g + {\a\over L}\right)\right)
\right|^2
= \left(1 - \left({\sin\pi\a\over\pi\a}\right)^2\right){T\over2\pi}
\log^2T + O(T\log T)
$$
holds uniformly for $|\a| \le \hf L$, where $L = {1\over2\pi}
\log({T\over2\pi})$. It would be interesting to recover this result
unconditionally, but our method of proof  does not
seem capable of achieving this.

\smallskip
If the RH holds, then $F(T,H) \equiv 0$ for $H>0$, and there is nothing more to say. However,
the RH is not known yet to hold, so that one may ask: what if RH fails, but $F(T,H) = 0\,$?
It follows that there exists a zeta zero $\b+i\g \;(T<\g\le T+H)$ such that $\b \ne \hf$.
If for such a zero one has $\hf < \b < 1$, then $1-\b+i\g$ is also a zero,
which follows from $\overline{\z(s)} = \z(\bar{s})$
and the functional equation
$$
\z(s) = \chi(s)\z(1-s), \quad \chi(s) := \frac{\G(\hf(1-s))}{\G(\hf s)}\pi^{s-1/2}\quad(\forall s\in \CC).
$$
Therefore one may consider only the case when $\hf < \b < 1$ and define, for a given $\g\;(>0)$,
$$
A(\g) := \sum_{\frac12<\b<1,\z(\b+i\g) = \z(\frac12+i\g)=0}1,\leqno(1.2)
$$
where the multiplicities of the zeros $\z(\b+i\g)$ are counted.
It is clear that
$$
0\;\le\; A(\g) \;\le\; N(\g + \hf) - N(\g-\hf) \;\ll\; \log\g.\leqno(1.3)
$$
It is reasonable to expect that $A(\g) =0$ for almost all $\g$, but this is not easy to prove.

As is customary, the function
$$
N(T) \= \sum_{0<\g\le T}1
$$
counts, with multiplicities, the number of zeta zeros, whose positive imaginary parts
do not exceed $T$. We have (see Chapter 1 of \cite{Iv1} or  Chapter 9 of \cite{Tit})
$$\eqalign{
N(T) &\,= \sum_{0<\g\le T}1 = {1\over\pi}\vartheta(T) + 1 + S(T),\cr
\vartheta(T) &:= \I\left\{\log\Gamma({\txt{1\over4}}+\hf iT)\right\}
- \hf T\log\pi,
\cr}\leqno(1.4)
$$
whence $\vartheta(T)$ is real and continuously differentiable. In fact, by
using Stirling's formula for the gamma-function, it is found that
$$
\vartheta(T) = {T\over2}\log{T\over2\pi} - {T\over2} - {\pi\over8}
+ O\left({1\over T}\right),\; \vartheta'(T) = \frac12 \log\frac{T}{2\pi} + O\Bigl(\frac{1}{T^2}\Bigr).
$$
Moreover,
$$
S(T) = {1\over\pi}\arg\z(\hf+iT) = {1\over\pi}\I\left\{
\log\z(\hf+iT)\right\} \ll \log T. \leqno(1.5)
$$
Thus (1.3) follows from (1.4) and (1.5).
Here for $T\not=\g$ the argument of $\z(\hf+iT)$ is
obtained by continuous variation along the straight lines joining the
points 2, $2+iT$, $\hf+iT$, starting with the value 0. If $T$ is  an
ordinate of a zeta-zero, then we define $S(T) = S(T+0)$. 
Clearly when $T \not= \g$ we can differentiate $S(T)$ by using (1.5). If $T=\g$,
then by (1.4) it is seen that $S(T)$ has a jump discontinuity which counts the
number of zeros $\rho$ with $\g = \I\rho = T$.
For a comprehensive account on $\z(s)$, the reader is referred to the monographs
 of E.C. Titchmarsh \cite{Tit} and the author \cite{Iv1}.

\medskip
There are some results for $F(T,T)$, defined by (1.1). The author \cite{Iv2} proved that unconditionally
$$
F(T,T) = \sum_{T<\g\le 2T}|\z(\hf+i\g)|^2 \,\ll_\e\, T\log^2T(\log\log T)^{3/2+\e},\leqno(1.6)
$$
where $\e$ denotes arbitrarily small positive numbers, not necessarily the same ones
at each occurrence, and $\ll_\e$ means that the
implied $\ll$-constant depends only on $\e$.
K. Ramachandra \cite{Ram} used a different method to obtain a result
which easily implies that the right-hand side of (1.6) is unconditionally $\ll T\log^2T\log\log T$.
The same bound was obtained by the author \cite{Iv3}, by another method. It was also used to
obtain several other results, namely
$$
\eqalign{&
\int_0^T|\zt|^2S(t)\d t \,\ll\, T\log T\log\log T,\cr&
\int_0^T|\zt|^2S^2(t)\d t \,\ll\, T\log T(\log\log T)^2,\cr}
\leqno(1.7)
$$
while under the Riemann Hypothesis one has
$$
\int_0^T|\zt|^2S(t)\d t \,\ll\, T\log T.\leqno(1.8)
$$

By a variant of the method used in \cite{Iv3} one can generalize these results to short intervals and
obtain the following unconditional results.

\medskip
THEOREM 1. {\it If $\,\,BT^{2/3}\log^4T  \le H = H(T) \le T\,$ for a suitable $B>0$, then we have}
$$
F(T,H) = \sum_{T<\g\le T+H}|\z(\hf+i\g)|^2 \,\ll\, H(\log T)^2\log \log T.
\leqno(1.9)
$$

\medskip
THEOREM 2. {\it If $\,\,BT^{2/3}\log^4T  \le H = H(T) \le T\,$ for a suitable $B>0$, then we have}
$$
\eqalign{&
\int_T^{T+H}|\zt|^2 S(t)\d t \,\ll\,H\log T\log\log T,\cr&
\int_T^{T+H}|\zt|^2 S^2(t)\d t \,\ll\,H\log T(\log\log T)^2.\cr}\leqno(1.10)
$$

\medskip
\head
2. The necessary lemmas
\endhead
If one defines
$$
R(t) := S(t) + \frac{1}{\pi}\sum_{p\le y}p^{-1/2}\sin(t\log p)\qquad(T\le t \le 2T), \leqno(2.1)
$$
where $p$ denotes primes, $y = T^\delta$, and $\delta >0$ is a small positive number, then
it is a classical result of A. Selberg \cite{Sel} that $R(t)$ is small on the average.
This was also later elaborated by K.-M. Tsang \cite{Tsa}. What is needed here is

\medskip
{\bf Lemma 1}. {\it Let $m>1$ be an integer, $1 < m\le (\log x)/192$, $x^{1/(4m)} < y \le x^{1/m}$
and $\,\log T \ll \log x \ll \log T$.
Then we have, for $\,T \ge T_0$,}
$$
\int_T^{T+H}R^{2m}(t)\d t < {(\E^{37}\pi^{-2}\e^{-3}m^2)}^m\,H\qquad(H = T^{27/82+\e}).   \leqno(2.2)
$$

\medskip
This is Lemma 7 of the paper of A.A. Karatsuba and M.A. Korolev \cite{KaKo}. Its good features are that
(2.2) is quite explicit, and moreover the range of $H$ is wide.

\medskip
{\bf Lemma 2}. {\it Let $\,BT^{2/3}\log^4T \le H \le T\,$ for a suitable $B>0$. Then}
$$
\eqalign{
\int_T^{T+H}|\zt|^4\d t &\ll H\log^4T,\cr
\int_T^{T+H}|\z'(\hf+it)|^4\d t & \ll H\log^8T.\cr}\leqno(2.3)
$$

\medskip
The first bound in (2.3) follows from the asymptotic formula
$$
\int_0^T|\zt|^4\d t = TP_4(\log T) + O(T^{2/3}\log^C T), \leqno(2.4)
$$
where $P_4(x)$ is a polynomial of degree four in $x$, with leading coefficient $1/(2\pi^2)$.
The proof of (2.4), with $C= 53/6$, was given by Y. Motohashi and the author \cite{IvMo}.
The value $C=8$ was given later by Y. Motohashi \cite{Mot2}.
The second bound in (2.3) follows from (4.1), (4.2) and (4.9) of the author's paper \cite{Iv4} and the first
bound in (2.3). It is clearly the range for $H$ in Lemma 2 which sets the limit to the range for
$H$ in Theorem 1 and Theorem 2.

\medskip
{\bf Lemma 3}. {\it Let $A(s) = \sum\limits_{m\le M}a(m)m^{-s}, \,a(m) \ll_\e m^\e$. Then
$$
\int_0^T|\zt A(\hf+it)|^2\d t = T\sum_{k,\ell\le M}\frac{a(k)\overline{a(\ell)}}{[k,\ell]}
\left(\log\frac{T(k,\ell)^2}{2\pi k\ell} + 2C_0 -1\right) + E(T,A),\leqno(2.5)
$$
where $C_0= -\G'(1)$ is Euler's constant, and $\,E(T,A) \ll_\e T^{1/3+\e}M^{4/3}$.}

\medskip
This is a version of the mean value theorem for a Dirichlet polynomial weighted by $|\zt|$,
and is due to Y. Motohashi \cite{Mot1}. As usual, $(k,\ell)$ is the greatest common divisor of
$k$ and $\ell$, and $[k,\ell]$ is their least common multiple.

\medskip
\head
3. Proof of Theorem 1
\endhead
Let henceforth $BT^{2/3}\log^4T \le H \le T$, and let $f(t)$ be a smooth function on $[T, T+H]$.
Then in view of (1.4) one has
$$
\eqalign{&
\sum_{T<\g\le T+H}f(\g) = \int_T^{T+H}f(t)\d N(t) \cr&=
\int_T^{T+H}f(t)\frac{1}{2\pi}\log \bigl(\frac{t}{2\pi}\bigr)\d t
+ \int_T^{T+H}f(t)\d \Bigl(S(t) + O\bigl(\frac1t\bigr)\Bigr) = I_1 + I_2,\cr}\leqno(3.1)
$$
say. The integral $I_1$ is usually not difficult to evaluate, and so is the
integral with $O(1/t)$, which is a continuously differentiable function. The main
problem is the evaluation of the integral in (3.1) with $S(t)$, which we write as
$$
I_2 = \int_T^{T+H}f(t)\d R(t) -
\int_T^{T+H}f(t)\frac{1}{\pi}\sum_{p\le T^\delta}p^{-1/2}\log p\cdot\cos(t\log p)\d t,\leqno(3.2)
$$
where (2.1) was used ($\delta>0$ is sufficiently small).
In the case of $f(t) \equiv |\zt|^2$, which is needed for Theorem 1, we easily see that
$$
I_1 \,\ll\, H(\log T)^2.\leqno(3.3)
$$
To deal with $I_2$, let
$$
{\Cal A}(T,H;V) := \Bigl\{\,t \,:\,(T\le t\le T+H)\,\wedge\,(|R(t)| \ge V)\,\Bigr\},
$$
where we suppose that $V = V(t) \ge 0$ and $\lim\limits_{T\to\infty}V(T) = +\infty$.
If $\mu(\cdot)$ denotes measure, then by Lemma 1 we obtain
$$
\eqalign{
\mu\Bigl({\Cal A}(T,H;V)\Bigr) & \,= \,\int_{{\Cal A}(T,H;V)}1\d t \,\le\, V^{-2m}\int_T^{T+H}R^{2m}(t)\d t
\cr& \ll \,\left(C\e^{-3}\Bigl(\frac{m}{V}\Bigr)^2\right)^m H,\cr}
$$
where $C, C_j,\ldots$ denote positive, absolute constants. If $m = [AV]$ for a sufficiently small
constant $A>0$, then
$$
\eqalign{
\left(C\e^{-3}\Bigl(\frac{m}{V}\Bigr)^2\right)^m & \le \left(C\e^{-3}A^2\right)^{[AV]}
\cr& \le\exp\left(-[AV]\log\frac{\e^3}{CA^2}\right)
\cr& \le\exp\left(-(AV -1)\log\frac{\e^3}{CA^2}\right) \le C_1\E^{-C_2V}\cr}
$$
for suitable $C_1, C_2$. If we choose
$$
V \= \frac{100}{C_2}\log\log T\qquad(T\ge T_0 > 0),\leqno(3.4)
$$
then we see that
$$
\mu\Bigl({\Cal A}(T,H;V)\Bigr) \ll H\E^{-C_2V} \= H(\log T)^{-100}.\leqno(3.5)
$$
Now we use (3.1) and (3.2) with $f(t) \equiv |\zt|^2 = \zt \z(\hf-it)$. This is needed since
integration by parts yields
$$
\eqalign{&
\int_T^{T+H}f(t)\d R(t)  \cr& = O(T^{1/3+\delta}) - \int_T^{T+H}R(t)\Bigl(\z'(\hf + it)\z(\hf-it)
- i\zt \z'(\hf-it) \Bigr)\d t.\cr}\leqno(3.6)
$$
Here we used the classical bound $\zt \ll |t|^{1/6}$ (see \cite{Iv1}).
Consider now the portion of the integral on the right-hand side of (3.6) for which $|R(t)| \ge V$,
where $V$ is given by (3.4).
By H\"older's inequality for integrals, this integral  does not exceed
$$
\eqalign{&
\left\{\mu\Bigl({\Cal A}(T,H;V)\Bigr)\int_T^{T+H}|\zt|^4\d t\int_T^{T+H}|\z'(\hf+it)|^4\d t
\int_T^{T+H}R^4(t)\d t\right\}^{1/4} \cr&\ll H,\cr}\leqno(3.7)
$$
on using (2.2) of Lemma 1 (with $m=2$), (2.3) of Lemma 2 and (3.5). The portion of the integral
over $[T, T+H] \,\backslash\, {\Cal A}(T,H;V)$ is
$$
\eqalign{&
\ll \log\log T \int_T^{T+H}|\zt||\z'(\hf+it)|\d t\cr&
\le \log\log T\left\{\int_T^{T+H}|\zt|^2\d t\int_T^{T+H}|\z'(\hf+it)|^2\d t\right\}^{1/2}\cr&
\ll \log\log T{\Bigl(H\log T\cdot H\log^3T\Bigr)}^{1/2} = H(\log T)^2\log\log T.\cr}
$$
The bounds for the mean square of $\z, \z'$ in short intervals follow similarly, but with less difficulty,
as the bound for the corresponding fourth moments in Lemma 2 (see e.g., Chapter 15 of \cite{Iv1}).

\medskip
It is also easily seen that (3.3) holds in our case. Thus it remains to estimate the second integral
in (3.2), namely
$$
I_3 := \frac{1}{\pi}\int_T^{T+H}|\zt|^2\sum_{p\le T^\delta}p^{-1/2}\log p\cdot\cos(t\log p)\d t.
\leqno(3.8)
$$
The integral in (3.8), by the Cauchy-Schwarz inequality, does not exceed
$$
\left\{\int_T^{T+H}|\zt|^2\d t\int_T^{T+H}|\zt|^2\Bigl|\sum_{p\le T^\delta}p^{-1/2-it}\log p\Bigr|^2\d t
\right\}^{1/2}.\leqno(3.9)
$$
As remarked above, the first integral in (3.9) is $\ll H\log T$, and for the second one we apply (2.5)
of Lemma 3, once with $T+H$ and once with $T$, and we subtract the results. In our interval for $H$ we shall
have $E(T+H, A) - E(T,A) \ll H$ for $M = T^\delta$ and sufficiently small $\delta>0$, where
$$
A(s) \;:=\; \sum_{p\le T^\delta}\log p\cdot p^{-1/2-it}.
$$
Further
$$
\eqalign{&
\sum_{p_1,p_2\le M}\frac{\log p_1\log p_2}{[p_1,p_2]}\left(\log \Bigl(\frac{T(p_1,p_2)^2}{2\pi p_1p_2}\Bigr)
+ 2C_0 -1\right)\cr&
= \sum_{p\le M}\frac{\log^2p}{p}\Bigl(\log \frac{T}{2\pi} + 2C_0-1\Bigr) \cr&
+ \sum_{p_1,p_2\le M; p_1\ne p_2}\frac{\log p_1\log p_2}{p_1p_2}
\Bigl(\log\Bigl(\frac{T}{2\pi p_1p_2}\Bigr)+2C_0-1\Bigr).
\cr}\leqno(3.10)
$$
The last expression is $\ll \log T\cdot\log^2M \ll \log^3T$,
if one uses the elementary bound
$$
\sum_{p\le x}\frac{\log p}{p} \ll \log x.
$$
Therefore the expression in (3.8) is
$\ll H\log^2T$, which finishes the proof of Theorem 1.

\medskip
\head
3. Proof of Theorem 2
\endhead
The proof of Theorem 2 is based on the same ideas as the proof of Theorem 1, so only its  salient
points will be mentioned. To prove the first bound in (1.10) we use (2.1). The integral with $R(t)$,
similarly as in the proof of Theorem 1, will be $\ll H\log T\log\log T$.  Let
$$
\sum(T) := \sum_{p\le T^\delta}p^{-1/2}\sin(t\log p).
$$
The contribution of $\sum(T)$ for which $|\sum(T)| \le \log\log T$ is trivially $$\ll H\log T\log\log T.$$
The remaining contribution is bounded by
$$
\frac{1}{\log\log T}\int_T^{T+H}|\zt|^2\Bigl|\sum_{p\le T^\delta}p^{-1/2-it}\Bigr|^2\d t.\leqno(4.1)
$$
The integral in (4.1) is estimated by Lemma 3, with the preceding $A(s)$ replaced by
$$
A_1(s) \;:=\; \sum_{p\le T^\delta} p^{-1/2-it}.
$$
This leads to an expression similar to the one in (3.9), namely
$$
\eqalign{&
\sum_{p_1,p_2\le M}\frac{1}{[p_1,p_2]}\left(\log \Bigl(\frac{T(p_1,p_2)^2}{2\pi p_1p_2}\Bigr)
+ 2C_0 -1\right)\cr&
= \sum_{p\le M}\frac{1}{p}\Bigl(\log \frac{T}{2\pi} + 2C_0-1\Bigr) \cr&
+ \sum_{p_1,p_2\le M; p_1\ne p_2}\frac{1}{p_1p_2}
\Bigl(\log\Bigl(\frac{T}{2\pi p_1p_2}\Bigr)+2C_0-1\Bigr) \ll \log T(\log\log T)^2,
\cr}\leqno(4.2)
$$
since $\sum\limits_{p\le x}1/p \ll \log\log x$. The bound in (4.2), combined with (4.1) leads then to
$$
\int_T^{T+H}|\zt|^2 S(t)\d t \,\ll\,H\log T\log\log T.
$$
To prove the remaining bound in (1.10) we use $S^2(t) \ll R^2(t) + \sum^2(t)$.
The integral with $R^2(t)$ is
$$
\ll H\log T(\log\log T)^2,\leqno(4.3)
$$
if we consider separately the cases $|R(t)| \le V$ and $|R(t)| \ge V$, where $V$ is as in (3.4).
The integral with $\sum^2(t)$ is estimated as the integral in (4.1). With the aid of (4.2) we
arrive again at the bound in (4.3), thereby completing the proof.

\vfill
\eject
\topglue1cm
\bigskip

\Refs
\medskip

\item{[Gon]} S.M. Gonek,  Mean values of the Riemann zeta-function and
its derivatives, Invent. math. {\bf75}(1984), 123-141.

\smallskip
\item{[Iv1]}
A. Ivi\'c, The Riemann zeta-function John Wiley \&
Sons, New York 1985 (reissue,  Dover, Mineola, New York, 2003).

\smallskip

\item{[Iv2]} A. Ivi\'c,  On certain sums over ordinates of zeta-zeros,
Bulletin CXXI de l'Acad\'emie Serbe des Sciences et des
Arts - 2001, Classe des Sciences math\'ematiques et naturelles,
Sciences math\'ematiques No. {\bf26}, pp. 39-52.

\smallskip
\item{[Iv3]} A. Ivi\'c,
On sums of squares of the Riemann zeta-function on the critical line,
Max-Planck-Institut f\"ur Mathematik, Preprint Series 2002({\bf52}), pp. 12.
(also in ``Bonner Math. Schriften'' Nr. {\bf360} (eds. D.R. Heath-Brown and B.Z. Moroz),
Proc. Session in analytic number theory and Diophantine equations, pp. 17).

\smallskip
\item{[Iv4]} A. Ivi\'c,
On certain moments of Hardy's function $Z(t)$ over short intervals, Moscow Journal
of Combinatorics and Number Theory, 2017, vol. 7, issue 2,
pp. 59-73, [pp. 147-161].

\smallskip
\item{[IvMo]} A. Ivi\'c and Y. Motohashi, On the fourth power moment of the
Riemann zeta-function, J. Number Theory {\bf51}(1995), 16-45.

\smallskip
\item{[KaKo]}
A.A. Karatsuba and M.A. Korolev,  The behavior of the argument of the
Riemann zeta function on the critical line. (Russian) Uspekhi Mat. Nauk {\bf61}(2006), no. 3(369), 3-92;
 translation in Russian Math. Surveys {\bf61}(2006), no. 3, 389-482.

\smallskip
\item{[Mot1]} Y. Motohashi, A note on the mean value of the zeta and $L$-functions,
Proc. Japan Acad. {\bf62}, Ser. {\bf A}(1986), 399-403.

\smallskip
\item{[Mot2]}
Y. Motohashi, Spectral  theory of the Riemann zeta-function,
Cambridge University Press, Cambridge, 1997.

\smallskip
\item{[Ram]} K. Ramachandra,  On a problem of Ivi\'c,
Hardy-Ramanujan Journal {\bf23}(2001), 10-19.

\smallskip
\item{[Sel]} A. Selberg, On the remainder formula for $N(T)$, the number of zeros of $\z(s)$
in the strip $0<t < T$, Avhandliger utgitt av Det Norske Videnskaps-Akademi i Oslo I.
Mat.-Naturv. Klasse (1944), No. 1, 1-27.

\smallskip
\item{[Tit]}  E.C. Titchmarsh,
The theory of the Riemann zeta-function, 2nd ed. edited by
D.R. Heath-Brown, Oxford, Clarendon Press, 1986.

\smallskip
\item{[Tsa]} K.-M. Tsang,  Some $\Omega$--theorems for the Riemann
zeta-function, Acta Arithmetica {\bf 46}(1986), 369-395.

\endRefs

\enddocument

\bye